\documentclass[12pt, reqno]{amsart}
\usepackage{amsmath, amsthm, amscd, amsfonts, amssymb, graphicx, color}
\usepackage[bookmarksnumbered, colorlinks, plainpages]{hyperref}

\makeatletter \oddsidemargin.9375in \evensidemargin \oddsidemargin
\def\authorsaddresses#1{\dedicatory{#1}}

\newtheorem{theorem}{Theorem}[section]
\newtheorem{lemma}[theorem]{Lemma}

\newtheorem{corollary}[theorem]{Corollary}
\theoremstyle{definition}
\newtheorem{definition}[theorem]{Definition}
\newtheorem{example}[theorem]{Example}

\theoremstyle{remark}
\newtheorem{remark}[theorem]{Remark}
\numberwithin{equation}{section}
\begin{document}

\title[Double controlled cone metric spaces]
{Double controlled cone metric spaces and the related fixed point theorems}
\author[Tayebeh Lal Shateri]{Tayebe Lal Shateri}
\authorsaddresses{ Tayebeh Lal Shateri\\ Department of Mathematics and Computer
Sciences,\\ Hakim Sabzevari University, Sabzevar, P.O. Box 397, IRAN\\
t.shateri@hsu.ac.ir; shateri@ualberta.ca}

\address{Tayebe Lal Shateri \\ Department of Mathematics and Computer
Sciences, Hakim Sabzevari University, Sabzevar, P.O. Box 397, IRAN}
\email{ \rm t.shateri@hsu.ac.ir; shateri@ualberta.ca}
\thanks{*The corresponding author:
t.shateri@hsu.ac.ir ; shateri@ualberta.ca (Tayebe Lal Shateri)}
 \maketitle

\begin{abstract}
In this paper, we introduce double controlled cone metric spaces via two control functions. An example of a double controlled cone metric space by two incomparable functions, which is not
a controlled metric space, is given. We also provide some fixed
point results involving Banach type and Kannan type contractions in the setting of double controlled cone metric spaces.
 \vskip 3mm
\end{abstract}

\section{Introduction and preliminaries}\vskip 2mm
One interesting extension of metric spaces is $b$-metric spaces introduced by Bakhtin \cite{BAK}. Recently, several extensions of $b$-metric spaces such as $b_v(s)$-metric spaces and $b$-rectangular metric spaces were introduced and some fixed point theorems are proved on these spaces, see (\cite{MI,RO}). In 2018, Mlaiki et al. \cite{ML1} gave an extension of the extended $b$-metric spaces, which was called controlled metric type spaces. Also, in 2018, Abdeljawad et al. \cite{ABD}, introduced the concept of double controlled metric type spaces. Also, they proved an analogue of the  Banach contraction principle on controlled metric type spaces. Banach contraction principle \cite{BA} plays an important role in several branches of mathematics. For instance, it has been used to study the existence of solutions for nonlinear Volterra integral equations and nonlinear integro-differential equations in Banach spaces and to prove the convergence of algorithms in computational mathematics. In 2020, Mlaiki \cite{ML2} gave an extension of all the types of metric spaces mentioned, which was called double controlled metric-like spaces, and proved some fixed point results on double controlled metric-like spaces. \\Fixed point theorems are the basic mathematical tools used in showing the existence of solution  in such diverse ﬁelds such as biology, chemistry, economics, engineering, and game theory.

Haung and Zhang \cite{HU} have introduced the concept of the cone metric space, replacing the set of real numbers by an ordered Banach space and they showed some fixed point theorems of contractive type mappings on cone metric spaces. Later many
authors generalized their results, see (\cite{AB1,AB2,CA,CAK,TU,XI}). 

In this paper, by using some ideas of \cite{ABD,HU}, we introduce double controlled cone metric spaces via two control functions. An example of a double controlled cone metric space by two incomparable functions, which is not a controlled metric space, is given. We also provide some fixed point results involving Banach type and Kannan type contractions in the setting of double controlled cone metric spaces.\\The following definitions and results will be needed in the sequel.

Let $E$ be a real Banach space and $P$ be a subset of $E$. $P$ is called a cone if it satisfies the followings\\
$(C_1)\; P \;\text{is closed, non-empty and}\; P\neq\{0\},$\\
$(C_2)\; ax+by\in P \quad \text{for all}\; x,y\in P\; \text{and non-negative real numbers}\; a,b,$\\
$(C_3)\; P\cap (-P)=\{0\}.$\\
For a given cone $P\subseteq E$, we can deﬁne a partial ordering $\preceq$ on $E$ with respect to $P$ by $x\preceq y$ if and only if $y-x\in P$. We write $x\prec y$ to 
indicate that $x\preceq y$ but $x\neq y$, while $x\ll y$ will stand for $y-x\in  intP$, in which $intP$ denotes the interior of $P$. 

Let $E$ be a Banach space, $P$ be a cone in $E$ with $intP\neq\Phi$ and $\preceq$ is partial ordering with respect to $P$. The cone $P$ is called normal if
\begin{equation}\label{eq0}
\inf\{\|x+y\|:x,y\in P,\; \|x\|=\|y\|=1\}\succ 0,
\end{equation}
or equivalently, there is a constant number $M\succ 0$ such that for all $x,y\in E$ where $0\preceq x\preceq y$ implies $\|x\|\preceq M\|y\|$. The least positive number satisfying above is called the normal constant of $P$. 
From \eqref{eq0} it is easy to see that $P$ is a non-normal if and only if there exist sequences $x_n,y_n\in P$ such that 
$$0\preceq x_n\preceq x_n+y_n\;,\quad x_n+y_n\rightarrow 0\quad \text{but} \quad x_n\nrightarrow 0.$$
We give an example of non-normal cone.
\begin{example}\cite{DE} 
Let $E=C^1[0,1]$ with the norm $\|x\|=\|x\|_{\infty}+\|x'\|_{\infty}$ on $P=\{x\in E:x(t)\geqslant 0 \;{\text on}\; [0,1]\}$. This cone is not normal. Consider
$$x_n(t)=\frac{1-\sin nt}{n+2},\; y_n(t)=\frac{1+\sin nt}{n+2}.$$
Since $\|x_n\|=\|y_n\|=1$ and $\|x_n+y_n\|=\frac{2}{n+2}\rightarrow 0$, it follows that $P$ is a non-normal cone.
\end{example}
\begin{definition}(Cone metric space) Let $\mathcal X$ be a non-empty set. The  mapping \\$d:\mathcal X\times\mathcal X\to E$ is said to be a cone metric on $\mathcal X$ if for all $x,y,z\in X$ the followings hold\\
$(CM_1)\; 0\prec d(x,y)\;\text{and}\; d(x,y)=0\; \text{if and only if}\; x=y,$\\
$(CM_2)\; d(x,y)=d(y,x),$\\
$(CM_3)\; d(x,y)\preceq d(x,z)+d(z,y).$\\
Then $(\mathcal X,d)$ is called a cone metric space.
\end{definition} 
We now state the definition of a controlled cone metric space.\begin{definition} (Controlled cone metric space) 
Let $\mathcal X$ be a nonempty set and let $\alpha:\mathcal X\times\mathcal X\to [1,\infty)$ be a function. $p:\mathcal X\times\mathcal X\to E$ is called a controlled cone metric type with $\alpha$ if the followings hold\\
$(CCM_1)\; 0\prec p(x,y) \;{\text and}\; p(x,y)=0\; \text {if and only if}\; x=y,$\\
$(CCM_2)\; p(x,y)=p(y,x),$\\
$(CCM_3)\; p(x,y)\preceq\alpha(x,z)p(x,z)+\alpha(z,y)p(z,y),$
\end{definition}
for all $x,y,z\in X$. Then the pair $(\mathcal X,p)$ is said to be a controlled cone metric space. Each cone metric space is a controlled cone metric space with $\alpha(x,y)=1$.
\begin{definition} (Double controlled cone metric space) 
Let $\mathcal X$ be a nonempty set and let $\alpha ,\beta:\mathcal X\times\mathcal X\to [1,\infty)$ be non-comparable functions. $p:\mathcal X\times\mathcal X\to E$ is called a double controlled cone metric with respect to $\alpha$ and $\beta$ if the followings hold\\
$(DCM_1)\; 0\prec p(x,y) \;{\text and}\; p(x,y)=0\; \text {if and only if}\; x=y,$\\
$(DCM_2)\; p(x,y)=p(y,x),$\\
$(DCM_3)\; p(x,y)\preceq\alpha(x,z)p(x,z)+\beta(z,y)p(z,y),$
for all $x,y,z\in X$. Then the pair $(\mathcal X,p)$ is said to be a double controlled cone metric space. 
\end{definition}
It is clear that each controlled cone metric space is a double controlled cone metric space, but there exists double controlled cone metric spaces which are not controlled cone metric space.
\begin{example}\label{ex1}
Let $E=\mathbb R^2$, $P=\{(x,y)\in E:\; x,y\geqslant 0\}$, $\mathcal X=\mathbb R^+$ and $p:\mathcal X\times \mathcal X\to E$ be deﬁned by 
\begin{equation*}
p(x,y) = \left\{
\begin{array}{rl}
    (0,0),  &   \Leftrightarrow x=y \\
    (\frac{1}{x},\frac{1}{3}),  &   x\geqslant 1, y\in [0,1)\\
     (\frac{1}{3},\frac{1}{y}),  &   x\in [0,1),y\geqslant 1\\
     (1,1), & otherwise.
\end{array} \right.
\end{equation*}
Consider $\alpha,\beta:\mathcal X\times\mathcal X\to [1,\infty)$ by
\begin{equation*}
\alpha(x,y) = \left\{
\begin{array}{rl}
    x,  &   x,y\geqslant 1 \\
    1,  &   otherwise
\end{array} \right. \quad {\text and}\qquad 
\beta(x,y) = \left\{
\begin{array}{rl}
   1,  &   x,y< 1 \\
   \max\{x,y\},  &   otherwise.
\end{array} \right.
\end{equation*}
Then $(\mathcal X,p)$ is a double controlled cone metric space. On the other hand, we have
$$p(0,\frac{1}{2})=(1,1)>(\frac{2}{3},\frac{2}{3})=\alpha(0,3)p(0,3)+\alpha(3,\frac{1}{2})p(3,\frac{1}{2}).$$
This implies that $p$ is not a controlled cone metric space when $\alpha=\beta$. 
\end{example}
\begin{example}\label{ex2}
Let $E=\mathbb R^2$, $P=\{(x,y)\in E:\; x,y\geqslant 0\}$ and let \\$\mathcal X=\{(x,0):0\leqslant x\leqslant 1\}\bigcup \{(0,x):0\leqslant x\leqslant 1\}$. Suppose that $p:\mathcal X\times \mathcal X\to E$ be deﬁned by
\begin{align*}
p\left((x,0),(y,0)\right)&=\left(\frac{4}{3}|x-y|,|x-y|\right),
p\left((0,x),(0,y)\right)=\left(|x-y|,\frac{2}{3}|x-y|\right)\\
p\left((x,0),(0,y)\right)&=p\left((0,y),(x,0)\right)=\left(\frac{4}{3}x+y,x+\frac{2}{3}y\right).
\end{align*}
Define $\alpha,\beta:\mathcal X\times\mathcal X\to [1,\infty)$ as
\begin{equation*}
\alpha(X,Y) = \left\{
\begin{array}{rl}
    \max\{\frac{1}{x},\frac{1}{y}\},  &   x,y\neq 0\\
    (1,1),  &   x=y=0
\end{array} \right. \quad {\text and}\qquad 
\beta(X,Y) = \left\{
\begin{array}{rl}
   \frac{1}{x}+\frac{1}{y},  &   x,y\neq 0 \\
   (1,1),  &   x=y=0,
\end{array} \right.
\end{equation*}
for all $X,Y\in \mathcal X$ in which $X=(x,0)$, $Y=(y,0)$ or $X=(0,x)$, $Y=(0,y)$. Then $p$ is a double controlled cone metric. In fact, the conditions $(DCM_1)$ and $(DCM_2)$ hold. $(DCM_3)$ is satisfied too, because by the example in \cite{HU}, $p$ is a cone metric space and so we have
$$p(x,y)\leqslant p(x,z)+p(y,z)\leqslant\alpha(x,z)p(x,z)+\beta(y,z)p(y,z).$$
\end{example}
\begin{definition}
Let $(\mathcal X,p)$ be a double controlled cone metric space with respect to $\alpha$ and $\beta$.\\ $(i)$ The sequence $\{x_n\}$ is convergent to some $x$ in 
$\mathcal X$, if for every $c\in E$ with $0\ll c$ there is $N$ such that for all $n>N$, $p(x_n,x)\ll c$, then $\{x_n\}$ is
said to be convergent and $\{x_n\}$ converges to $x$, and $x$ is the limit of $\{x_n\}$. It is written as $\lim_{n\to\infty}x_n=x.$\\
$(ii)$ The sequence $\{x_n\}$ is said Cauchy, if for every $c\in E$ with $0\ll c$ there is $N$ such that for all $m,n>N$, $p(x_m,x_n)\ll c$.\\
$(iii)
\; (\mathcal X,p)$ is said complete if every Cauchy sequence is convergent.
\end{definition}
The following result is similar to \cite[Lemma 1]{HU}, and so the proof is removed.
\begin{lemma}
Let $(\mathcal X,p)$ be a double controlled cone metric space with respect to $\alpha$ and $\beta$, $P$ be a normal cone with normal constant $M$. Let $\{x_n\}$ be a sequence in $\mathcal X$. Then $\{x_n\}$ converges to $x$ if and only if $\lim_{n\to\infty}p(x_n,x)=0$.
\end{lemma}
\begin{lemma}
Let $(\mathcal X,p)$ be a double controlled cone metric space with respect to $\alpha$ and $\beta$, $P$ be a normal cone with normal constant $M$. Let $\{x_n\}$ be a sequence in $\mathcal X$ such that $\{x_n\}$ converges to $x$ and $y$. If $\lim_{n\to\infty}\alpha(x,x_n)$ and $\lim_{n\to\infty}\beta(x_n,y)$ exist and are finite, then $x=y$. 
That is the limit of $\{x_n\}$ is unique.
\begin{proof}
For any $c\in E$ with $0\ll c$, there is $N$ such that for all $n>N$, $p(x_n,x)\ll c$ and $p(x_n,y)\ll c$. We have
$$p(x,y)\preceq \alpha(x,x_n)p(x,x_n)+\beta(x_n,y)p(x_n,y)\preceq c\big(\alpha(x,x_n)+\beta(x_n,y)\big).$$
Hence, $\|p(x,y)\|\preceq c\|\alpha(x,x_n)+\beta(x_n,y)\|$. Since $\lim_{n\to\infty}\alpha(x,x_n)$ and $\lim_{n\to\infty}\beta(x_n,y)$ exist and are finite and $c$ is arbitrary, so $p(x,y)=0$,  therefore $x=y$.
\end{proof}
\end{lemma}
\section{ Main results}
In this section, we give some fixed point results in double controlled cone metric spaces. We assume that $(\mathcal X,p)$ be a complete double controlled cone metric space with respect to the functions $\alpha,\beta: \mathcal X\times \mathcal X\to[1,\infty)$ and $P$ be a normal cone with normal constant $M$.
\begin{theorem}\label{th1}
Let $T:\mathcal X\to\mathcal X$ satisfies the contraction condition\begin{equation}\label{eq1.1}
p(Tx,Ty)\preceq kp(x,y),
\end{equation}
for all $x, y\in\mathcal X$, where $k\in(0,1)$. For $x_0\in\mathcal  X$, choose $x_n=T^nx_0$. Suppose that
\begin{equation}\label{eq1.2}
\sup_{m\succeq 1}\lim_{i\to\infty}\frac{\alpha(x_{i+1},x_{i+2})}{\alpha(x_i,x_{i+1})}\beta(x_{i+1},x_m)\prec\frac{1}{k}.
\end{equation}
If for each $x\in\mathcal X,
\;\lim_{n\to\infty}\alpha(x,x_n)$ and $\lim_{n\to\infty}\beta(x_n,x)$ exist and are finite, then $T$ has a unique fixed point.
\begin{proof}
Suppose that the sequence $\{x_n\}$ in $\mathcal X$ satisfies the hypothesis of the theorem. Then by \eqref{eq1.1}, we get
$$p(x_n,x_{n+1})\preceq k^np(x_1,x_0),$$
for all $n\succeq 0$. Let $m,n$ be integers such that $m\succ n$.  Similar to proof of \cite[Theorem 1]{ABD}, we have
\begin{align*}
p(x_n,x_m)&\preceq\alpha(x_n,x_{n+1})p(x_n,x_{n+1})+\beta(x_{n+1},x_m)p(x_{n+1},x_m)\\
&\preceq\alpha(x_n,x_{n+1})p(x_n,x_{n+1})+\beta(x_{n+1},x_m)\alpha(x_{n+1},x_{n+2})p(x_{n+1},x_{n+2})\\
&+\beta(x_{n+1},x_m)\beta(x_{n+2},x_m)p(x_{n+2},x_m)\\
&\preceq\cdots\\
&\preceq\alpha(x_n,x_{n+1})p(x_n,x_{n+1})+\sum_{i=n+1}^{m-2}\Big(\prod_{j=n+1}^{i}\beta(x_j,x_m)\Big)\alpha(x_i,x_{i+1})p(x_i,x_{i+1})\\&+\prod_{k=n+1}^{m-1}\beta(x_k,x_m)p(x_{m-1},x_m)\\
&\preceq\alpha(x_n,x_{n+1})k^np(x_1,x_0)+\sum_{i=n+1}^{m-2}\Big(\prod_{j=n+1}^{i}\beta(x_j,x_m)\Big)\alpha(x_i,x_{i+1})k^ip(x_1,x_0)\\&+\prod_{k=n+1}^{m-1}\beta(x_k,x_m)k^{m-1}p(x_1,x_0)\\
&\preceq\alpha(x_n,x_{n+1})k^np(x_1,x_0)+\sum_{i=n+1}^{m-1}\Big(\prod_{j=n+1}^{i}\beta(x_j,x_m)\Big)\alpha(x_i,x_{i+1})k^ip(x_1,x_0)\\
&\preceq\alpha(x_n,x_{n+1})k^np(x_1,x_0)+\sum_{i=n+1}^{m-1}\Big(\prod_{j=0}^{i}\beta(x_j,x_m)\Big)\alpha(x_i,x_{i+1})k^ip(x_1,x_0),
\end{align*}
hence
\begin{align}\label{eq1.3}
\|p(x_n,x_m)\|&\preceq M\|\alpha(x_n,x_{n+1})k^np(x_0,x_1)\nonumber\\&+\sum_{i=n+1}^{m-1}\Big(\prod_{j=0}^{i}\beta(x_j,x_m)\Big)\alpha(x_i,x_{i+1})k^ip(x_0,x_1)\|.
\end{align}
Set $S_r=\sum_{i=0}^{r}\Big(\prod_{j=0}^{i}\beta(x_j,x_m)\Big)\alpha(x_i,x_{i+1})k^i$, then we have $$\|p(x_n,x_m)\|\preceq M\|p(x_0,x_1)\big[k^n\alpha(x_n,x_{n+1})+(S_{m-1}-S_n)\big]\|.$$
\eqref{eq1.2} implies that the limit of the real sequence $\{S_n\}$ exists and so $\{S_n\}$ is Cauchy. Letting $m,n\to\infty$ in \eqref{eq1.3}, yields $\lim_{m,n\to\infty}p(x_n,x_m)=0$, and so
$\{x_n\}$ is a Cauchy sequence. By the completeness of $\mathcal X$, there exists $x\in\mathcal X$ such that $\lim_{n\to\infty}x_n=x$. We claim that $Tx=x$. It follows from $(DCM_3)$ and \eqref{eq1.1} that
\begin{align*}
p(x,Tx)&\preceq\alpha(x,x_{n+1})p(x,x_{n+1})+\beta(x_{n+1},Tx)p(x_{n+1},Tx)\\
&\preceq\alpha(x,x_{n+1})p(x,x_{n+1})+k\beta(x_{n+1},Tx)p(x_n,x),
\end{align*}
and so $$\|p(x,Tx)\|\preceq M\Big(\|\alpha(x,x_{n+1})p(x,x_{n+1})\|+k\|\beta(x_{n+1},Tx)p(x_n,x)\|\Big).$$
Hence, $\|p(x,Tx)\|=0$, that is, $Tx=x$. Now, if $y$ is another fixed point of $T$, then
$$p(x,y)=p(Tx,Ty)\preceq kp(x,y),$$
thus $\|p(x,y)\|=0$, and so $x=y$. Therefore the fixed point of $T$ is unique.
\end{proof}
\end{theorem}
Now, we give examples of mappings on a double controlled cone metric space which satisfy in the conditions of Theorem \ref{th1}.
\begin{example}
Let $E=\mathbb R^2$, $P=\{(x,y)\in E:\; x,y\geqslant 0\}$ and let \\$\mathcal X=\{(x,0):0\leqslant x\leqslant 1\}\bigcup \{(0,x):0\leqslant x\leqslant 1\}$. Suppose that $p:\mathcal X\times \mathcal X\to E$ be deﬁned by
\begin{align*}
p\left((x,0),(y,0)\right)&=\left(\frac{4}{3}|x-y|,|x-y|\right),
p\left((0,x),(0,y)\right)=\left(|x-y|,\frac{2}{3}|x-y|\right)\\
p\left((x,0),(0,y)\right)&=p\left((0,y),(x,0)\right)=\left(\frac{4}{3}x+y,x+\frac{2}{3}y\right).
\end{align*}
Define $\alpha,\beta:\mathcal X\times\mathcal X\to [1,\infty)$ as
$$\alpha(X,Y)=\beta(X,Y)=1,$$
for all $X,Y\in\mathcal X$. Then $p$ is a double controlled cone metric. Define $T:\mathcal X\to\mathcal X$ be defined by
$$T(x,0)=(\frac{x}{2},0),\; T(0,y)=(0,\frac{y}{2}),$$
for all $x,y\in[0,1]$. Then $T$ satisfies condition \eqref{eq1.1}. Note that for each $X_0\in\mathcal X$, $X_0=(x_0,0)$ or $X_0=(0,x_0)$, and then  $X_n=T^nX_0=(\frac{x_0}{2^n},0)$ or $X_n=T^nX_0=(0,\frac{x_0}{2^n})$. All hypotheses of Theorem \ref{th1} are satisfied with $k=\frac{1}{2}$. In fact, $(0,0)$ is the unique fixed point of $T$.
\end{example}
The next result has an analogue on metric spaces \cite{KAR1}.
\begin{theorem}\label{th2}
Suppose that the mapping $T:\mathcal X\to\mathcal X$ satisfies the contractive condition
\begin{equation}\label{eq2.1}
p(Tx,Ty)\preceq ap(x,Tx)+bp(y,Ty)
\end{equation}
for all $x,y\in\mathcal X$, where $a+b\prec 1\; {\text and}\;a,b\in[0,1)$. For arbitray $x_0\in\mathcal  X$, choose $x_n=T^nx_0$. Assume that
\begin{equation}\label{eq2.2}
\sup_{m\succeq 1}\lim_{i\to\infty}\frac{\alpha(x_{i+1},x_{i+2})}{\alpha(x_i,x_{i+1})}\beta(x_{i+1},x_m)\prec\frac{1-b}{a}.
\end{equation}
If for each $x\in\mathcal X$,
\begin{equation}\label{eq2.3}
\lim_{n\to\infty}\alpha(x,x_n)\; \text{exists},\; \text{finite and}\quad \lim_{n\to\infty}\beta(x_n,x)\prec\frac{1}{b},
\end{equation}
then $T$ has a unique fixed point in X.
\begin{proof}
Let the sequence $\{x_n\}$ in $\mathcal X$ satisfy the hypothesis of the theorem. From \eqref{eq2.1}, for all $n\succeq 0$ we obtain
\begin{align*}
p(x_n,x_{n+1})&=p(Tx_{n-1},Tx_n)\\
&\preceq ap(x_{n-1},Tx_{n-1})+bp(x_n,Tx_n)\\
&=ap(x_{n-1},xn)+bp(x_n,x_{n+1}),
\end{align*}
hence $p(x_n,x_{n+1})\preceq\frac{a}{1-b}p(x_{n-1},x_n)$. Continuing this process we get
\begin{equation}\label{eq2.4}
p(x_n,x_{n+1})\preceq\Big(\frac{a}{1-b}\Big)^np(x_1,x_0),\quad \text{for all}\; n\succeq 0.
\end{equation}
Now, we show that the sequence $\{x_n\}$ is Cauchy. Using $(DCM_3)$ and \eqref{eq2.4}, for all $m,n\in\mathbb N$ we get
\begin{align*}
p(x_n,x_m)&\preceq\alpha(x_n,x_{n+1})p(x_n,x_{n+1})+\beta(x_{n+1},x_m)p(x_{n+1},x_m)\\
&\preceq\cdots\\
&\preceq\alpha(x_n,x_{n+1})p(x_n,x_{n+1})+\sum_{i=n+1}^{m-2}\Big(\prod_{j=n+1}^{i}\beta(x_j,x_m)\Big)\alpha(x_i,x_{i+1})p(x_i,x_{i+1})\\&+\prod_{k=n+1}^{m-1}\beta(x_k,x_m)p(x_{m-1},x_m)\\
&\preceq\alpha(x_n,x_{n+1})\Big(\frac{a}{1-b}\Big)^np(x_1,x_0)+\sum_{i=n+1}^{m-2}\Big(\prod_{j=n+1}^{i}\beta(x_j,x_m)\Big)\alpha(x_i,x_{i+1})\Big(\frac{a}{1-b}\Big)^ip(x_1,x_0)\\&+\prod_{k=n+1}^{m-1}\beta(x_k,x_m)\Big(\frac{a}{1-b}\Big)^{m-1}p(x_1,x_0).
\end{align*}
Since $a+b\prec 1$, so $\frac{a}{1-b}\prec 1$. This implies that $\|p(x_n,x_m)\|\to 0$ as $m,n\to\infty$, therefore the sequence $\{x_n\}$ is Cauchy, and the completeness of $\mathcal X$ implies that there exists an element $x\in\mathcal X$ such that $\{x_n\}$  converges to $x$. If $Tx\neq x$, we deduce that
\begin{align*}
0\prec p(x,Tx)&\preceq\alpha(x,x_{n+1})p(x,x_{n+1})+\beta(x_{n+1},Tx)p(x_{n+1},Tx)\\
&\preceq\alpha(x,x_{n+1})p(x,x_{n+1})+\beta(x_{n+1},Tx)\big[ap(x_n,x_{n+1})+bp(x,Tx)\big],
\end{align*}
hence
\begin{equation}\label{eq2.5}
0\prec\|p(x,Tx)\|\preceq M\|\alpha(x,x_{n+1})p(x,x_{n+1})+\beta(x_{n+1},Tx)\big[ap(x_n,x_{n+1})+bp(x,Tx)\big]\|.
\end{equation}
Making use of the condition \eqref{eq2.3}, and passing to the limit on \eqref{eq2.5}, we get $$0\prec\|p(x,Tx)\|\prec\|p(x,Tx)\|$$ which is a contradiction, therefore $Tx=x$. Suppose that $T$ has another fixed point $y$, then
\begin{align*}
p(x,y)&=p(Tx,Ty)\preceq ap(x,Tx)+bp(y,Ty)\\
&=ap(x,x)+bp(y,y)=0.
\end{align*}
Consequently $x=y$, and $T$ has a unique fixed point.
\end{proof}
\end{theorem}
If in Theorem \ref{th2}, we consider $a=b\in[0,\frac{1}{2})$, then we get the following result.
\begin{corollary}
Suppose that the mapping $T:\mathcal X\to\mathcal X$ satisfies the contractive condition
\begin{equation}
p(Tx,Ty)\preceq a\big[p(x,Tx)+p(y,Ty)\big]
\end{equation}
for all $x,y\in\mathcal X$, where $a+b\prec 1\; {\text an}\;a,b\in[0,1)$. For arbitrary $x_0\in\mathcal  X$, choose $x_n=T^nx_0$. Assume that
\begin{equation}
\sup_{m\succeq 1}\lim_{i\to\infty}\frac{\alpha(x_{i+1},x_{i+2})}{\alpha(x_i,x_{i+1})}\beta(x_{i+1},x_m)\prec\frac{1-a}{a}.
\end{equation}
If for each $x\in\mathcal X$,
\begin{equation}
\lim_{n\to\infty}\alpha(x,x_n)\; \text{exists},\; \text{is finite and}\quad \lim_{n\to\infty}\beta(x_n,x)\prec\frac{1}{a},
\end{equation}
then $T$ has a unique fixed point in X.
\end{corollary}
Following is the Reich \cite{RE} type contraction mapping considered here to prove  another fixed point theorem in double controlled cone metric spaces.
\begin{theorem}\label{th3}
Let the mapping $T:\mathcal X\to\mathcal X$ satisfy the contractive
condition 
\begin{equation}\label{eq3.1}
p(Tx,Ty)\preceq ap(x,Tx)+bp(y,Ty)+cp(x,y)
\end{equation}
for all $x,y\in\mathcal X$, where $a+b+c\prec 1\;{\text and}\; a,b,c\in[0,1)$. For $x_0\in\mathcal  X$, choose $x_n=T^nx_0$. Assume that
\begin{equation}\label{eq3.2}
\sup_{m\succeq 1}\lim_{i\to\infty}\frac{\alpha(x_{i+1},x_{i+2})}{\alpha(x_i,x_{i+1})}\beta(x_{i+1},x_m)\prec\frac{1-b}{a+c}.
\end{equation}
If for each $x\in\mathcal X$,
\begin{equation}\label{eq3.3}
\lim_{n\to\infty}\alpha(x,x_n)\; \text{exists},\; \text{is finite and}\quad \lim_{n\to\infty}\beta(x,x_n)\prec\frac{1}{b},
\end{equation}
then $T$ has a unique fixed point in X.
\begin{proof}
Let the sequence $\{x_n\}$ in $\mathcal X$ satisfy the hypothesis of the theorem. From \eqref{eq3.1}, for all $n\succeq 0$ we obtain
\begin{align*}
p(x_n,x_{n+1})&=p(Tx_{n-1},Tx_n)\\
&\preceq ap(x_{n-1},Tx_{n-1})+bp(x_n,Tx_n)+cp(x_{n-1},x_n)\\
&=ap(x_{n-1},x_n)+bp(x_n,x_{n+1})+cp(x_{n-1},x_n),
\end{align*}
hence $p(x_n,x_{n+1})\preceq\frac{a+c}{1-b}p(x_{n-1},x_n)$. By induction we get
\begin{equation}\label{eq3.4}
p(x_n,x_{n+1})\preceq\Big(\frac{a+c}{1-b}\Big)^np(x_1,x_0),\quad \text{for all}\; n\succeq 0.
\end{equation}
Similar to the proof of Theorem \ref{th2}, we can show that the sequence $\{x_n\}$ is Cauchy. Infact $(DCM_3)$ and \eqref{eq3.4}, for all $m,n\in\mathbb N$ implies that
\begin{align*}
p(x_n,x_m)&\preceq\alpha(x_n,x_{n+1})\Big(\frac{a+c}{1-b}\Big)^np(x_1,x_0)\\&+\sum_{i=n+1}^{m-2}\Big(\prod_{j=n+1}^{i}\beta(x_j,x_m)\Big)\alpha(x_i,x_{i+1})\Big(\frac{a+c}{1-b}\Big)^ip(x_1,x_0)\\&+\prod_{k=n+1}^{m-1}\beta(x_k,x_m)\Big(\frac{a+c}{1-b}\Big)^{m-1}p(x_1,x_0).
\end{align*}
Since $a+b+c\prec 1$, so $\frac{a+c}{1-b}\prec 1$, this implies that $\|p(x_n,x_m)\|\to 0$ as $m,n\to\infty$. Therefore the sequence $\{x_n\}$ is Cauchy, and the completeness of $\mathcal X$ implies that there exists an element $x\in\mathcal X$ such that $\{x_n\}$  converges to $x$. If $Tx\neq x$, we deduce that
\begin{align*}
0\prec p(x,Tx)&\preceq\alpha(x,x_{n+1})p(x,x_{n+1})+\beta(x_{n+1},Tx)p(x_{n+1},Tx)\\
&\preceq\alpha(x,x_{n+1})p(x,x_{n+1})+\beta(x_{n+1},Tx)\big[ap(x_n,x_{n+1})+bp(x,Tx)+cp(x,x_n)\big],
\end{align*}
hence
\begin{equation}\label{eq3.5}
\|p(x,Tx)\|\preceq M\|\alpha(x,x_{n+1})p(x,x_{n+1})+\beta(x_{n+1},Tx)\big[ap(x_n,x_{n+1})+bp(x,Tx)+cp(x,x_n)\big]\|.
\end{equation}
Making use of the condition \eqref{eq3.3}, and passing to the limit on \eqref{eq3.5}, we get $$0\prec\|p(x,Tx)\|\prec\|p(x,Tx)\|$$ which is a contradiction, therefore $Tx=x$. Suppose that $T$ has another fixed point $y$, then
\begin{align*}
p(x,y)&=p(Tx,Ty)\preceq ap(x,Tx)+bp(y,Ty)+cp(x,y)\\
&=ap(x,x)+bp(y,y)+cp(x,y)=cp(x,y),
\end{align*}
since $c\prec 1$, hence $x=y$, and so $T$ has a unique fixed point.
\end{proof}
\end{theorem}
If in Theorem \ref{th3}, we consider $a=b=c\in[0,\frac{1}{3})$, then we get the following result.
\begin{corollary}
Suppose that the mapping $T:\mathcal X\to\mathcal X$ satisfies the contractive condition
\begin{equation}
p(Tx,Ty)\preceq a\big[p(x,Tx)+p(y,Ty)+p(x,y)\big]
\end{equation}
for all $x,y\in\mathcal X$, where $a\in[0,\frac{1}{3})$. For arbitrary $x_0\in\mathcal  X$, choose $x_n=T^nx_0$. Assume that
\begin{equation}
\sup_{m\succeq 1}\lim_{i\to\infty}\frac{\alpha(x_{i+1},x_{i+2})}{\alpha(x_i,x_{i+1})}\beta(x_{i+1},x_m)\prec\frac{1-a}{2a}.
\end{equation}
If for each $x\in\mathcal X$,
\begin{equation}
\lim_{n\to\infty}\alpha(x,x_n)\; \text{exists},\; \text{is finite and}\quad \lim_{n\to\infty}\beta(x_n,x)\prec\frac{1}{a},
\end{equation}
then $T$ has a unique fixed point in X.
\end{corollary}
\begin{remark}
All the results which proved in this section, are hold in complete controlled cone metric spaces, too.
\end{remark}
\begin{remark}
In the case when the cone is not necessarily normal, 
we can prove Theorem \ref{th1}. For the proof we give the following properties of cone metrics which are often useful, when the cone is not normal.\\
$(i)$ If $x\preceq y$ and $y\ll z$, then $x\ll z$.\\
$(ii)$ If $c\in intP$, $0\preceq x_n$ and $x_n\rightarrow 0$, then there exists $N$ such that, for all $n>N$, we have $x_n\ll c$.\\
Now, we proof Theorem \ref{th1}, on a non-normal cone.
\begin{proof}
By the similar proof of Theorem \ref{th1} we have
$$p(x_n,x_m)\preceq \alpha(x_n,x_{n+1})k^np(x_0,x_1)+\sum_{i=n+1}^{m-1}\Big(\prod_{j=0}^{i}\beta(x_j,x_m)\Big)\alpha(x_i,x_{i+1})k^ip(x_0,x_1),$$
for $m\succ n$. Set $S_r=\sum_{i=0}^{r}\Big(\prod_{j=0}^{i}\beta(x_j,x_m)\Big)\alpha(x_i,x_{i+1})k^i$. Then we have $$p(x_n,x_m)\preceq p(x_0,x_1)\big[k^n\alpha(x_n,x_{n+1})+(S_{m-1}-S_n)\big].$$
\eqref{eq1.2} implies that the limit of the real sequence $\{S_n\}$ exists and so $\{S_n\}$ is Cauchy, hence $p(x_0,x_1)\big[k^n\alpha(x_n,x_{n+1})+(S_{m-1}-S_n)\big]\to 0$ as $m,n\to\infty$. From $(i)$ and $(ii)$ we deduce that $\{x_n\}$ is a Cauchy sequence. By the completeness of $\mathcal X$, there exists $x\in\mathcal X$ such that $\lim_{n\to\infty}x_n=x$. We claim that $Tx=x$. Since $T$ is continuous, so $\lim_{n\to\infty}Tx_n=Tx$. It follows from $(DCM_3)$, \eqref{eq1.1} and \eqref{eq1.3} that
\begin{align*}
p(x,Tx)&\preceq\alpha(x,x_n)p(x,x_n)+\beta(x_n,Tx)p(x_n,Tx)\\
&\preceq\alpha(x,x_n)p(x,x_n)+k\beta(x_n,Tx)p(x_{n-1},x).
\end{align*}
Given $c\in intP$, since $\lim_{n\to\infty}x_n=x$, so there exists $n\in\mathbb N$ such that, for all $n>N$, $p(x,x_{n-1})\ll\frac{c}{k}$ and $p(x_n,x)\ll c$. Hence, we obtain
\begin{align*}
p(x,Tx)\ll\alpha(x,x_n)c+k\beta(x_n,Tx)\frac{c}{k}=c\big[\alpha(x,x_n)+\beta(x_n,Tx)\big].
\end{align*}
Since $c$ is arbitrary, \eqref{eq1.3} implies that $p(x,Tx)=0$, and so $x$ is the fixed point of $T$.
\end{proof}
\end{remark}
\subsection*{Compliance with Ethical Standards:}
This article does not contain any studies with human participants or animals performed by author.

\end{document}